\documentclass[12pt]{article}
\usepackage[T2A]{fontenc}
\usepackage[cp1251]{inputenc}
\usepackage[english]{babel}
\usepackage{amssymb,amsmath}
\begin{document}

\begin{center}
\textsc{The Nilpotency Criterion for the }

\textsc{Derived Subgroup of a Finite Group}
\end{center}

\begin{center}
\textbf{Victor~S.~Monakhov}
\end{center}

\bigskip

{\bf Abstract.}
It is proved that the derived subgroup of a finite group is nilpotent
if and only if $|ab|\ge |a||b|$ for all primary commutators
$a$ and~$b$ of coprime orders.

\medskip


\bigskip

Bastos and Shumyatsky~\cite{BS} received the following
sufficient condition for nilpotency of the derived subgroup.

\medskip

{\bf Theorem~1.} {\sl Let $G$ be a finite group in which
$|ab|=|a||b|$ whenever the elements $a$ and $b$ have coprime orders.
Then $G$ is nilpotent.}

\medskip

Here the symbol $|x|$ denotes the order of the element~$x$ in a group~$G$,
$G^\prime$ is the derived subgroup of $G$.
In Theorem~1 we can limit to only primary commutators,
and we also obtain the nilpotency  criterion for the derived subgroup.
An element is said to be primary if it is of prime power order.

The aim of the present paper is to prove the following theorem.

\medskip

{\bf Theorem 2.} {\sl The derived subgroup of a finite group is nilpotent
if and only if $|ab|\ge |a||b|$ when  $a$ and~$b$ are the primary
commutators of coprime orders.}

\medskip

By $N_G(H)$ and~$C_G(H)$ we denote the normalizer  and centralizer
of a subgroup~$H$ in a group~$G$, respectively. If a group $G$ has
a normal Sylow $p$-subgroup, then $G$ is $p$\nobreakdash-\hspace{0pt}closed;
a group $G$ is $p$-nilpotent if there is a normal complement for
a Sylow $p$-subgroup in $G$. By $H\leftthreetimes K$ we denote
the semidirect product of a normal subgroup $H$ and a subgroup $K$.

\medskip

{\bf Lemma 3.} {\sl
Suppose that in a finite group $G$ the product of any two
primary commutators of coprime orders $s$ and $t$ is of order $\ge st$.
Let $H$ be a primary subgroup, $x$ be a primary commutator,
${(|x|,|H|)=1}$. If $x\in N_G(H)$, then $x\in C_G(H)$.}

\medskip

{\sc Proof.}
Let $y$ be an element of $H$. Since $[x,y]\in H$, the order of $[x,y]$ is
relatively prime to the order of $x$. Hence $|x[x,y]|\ge |x||[x,y]|$. As
$x[x,y]=y^{-1}xy$, we have $|x[x,y]|=|x|$. Now $[x,y]=1$ and $x\in C_G(H)$.
$\hfill \square$

\medskip

By $X_p(G)$ we denote the set of all commutators $x=[g,h]$ of~$G$
such that $|x|=p^n$ for a prime $p$ and a non-negative integer $n$.

\medskip

{\bf Proof of Theorem~2.}
Evidently, if the derived subgroup of a finite group is nilpotent, then
$|ab|=|a||b|$ for all primary commutators~$a$ and~$b$ of coprime orders.

Conversely, suppose that $G$ is a finite group such that $|ab|\ge |a||b|$
if elements $a$ and~$b$ are the primary commutators of coprime orders.
Firstly, we prove that $G$ is soluble.
Assume the converse and let $G$ be a counterexample of least order.
Then all proper subgroups of $G$ are soluble,
$G=G^\prime$ and $G/\Phi (G)$ is a simple group.
Here $\Phi (G)$ is the Frattini subgroup of~$G$.
Let $P$ be a Sylow $p$-subgroup of $\Phi (G)$,
$Q$ be a Sylow $q$-subgroup of~$G$, $p\ne q$.
By the Focal Subgroup Theorem \cite[7.3.4]{Gor},
$Q$ is generated by elements $x^{-1}y$, where $x,y\in Q$
and $x$ conjugates to $y$ in~$G$. Therefore
$$
x=y^g, \ \ g\in G, \ \
x^{-1}y=g^{-1}y^{-1}gy=[g,y]\in X_q(G),   \eqno (1)
$$
i.\,e. $Q$  is generated by elements of~$Q\cap X_q(G)$.
If $x\in Q\cap X_q(G)$, then $x\in C_G(P)$ in view of Lemma~3,
and $Q\le C_G(P)$. This is true for any $q\ne p$,
therefore $|G:C_G(P)|=p^a$. Since $C_G(P)$ is normal in $G$,
we have $G=C_G(P)$. This is true for any Sylow subgroup of~$\Phi (G)$,
so $\Phi (G)=Z(G)$.

Since $G$ is not $2$-nilpotent, it follows that in $G$
there is a $2$-closed Schmidt subgroup
(a nonnilpotent group all of whose proper subgroups
are nilpotent) $S=P\leftthreetimes Q$ of even order~\cite[IV.5.4]{Hup},
where $P$ is a normal in~$S$ Sylow $2$-subgroup,
$Q=\langle y\rangle$ is a nonnormal in $S$ cyclic Sylow $q$-subgroup.
The structure of $P$ is well known, in particular,
$Z(P)=P^\prime =\Phi (P)$ is an elementary abelian $2$-group.
Hence the exponent of~$P$ is at most~4.
Furthermore, $Z(S)=\Phi (P)\times \langle y^q\rangle$,
$S^\prime =P$~\cite[III.5.2]{Hup}. The last equality implies
that $X_2(S)$ is not contained in $Z(G)$.
If  $\forall a,b,g\in P$, then $[a,g][b,g]=[ab,g]$,
as $P^\prime = Z(P)$. Hence $P^\prime \subseteq X_2(S)$.
If~$P^\prime\nsubseteq Z(G)$, then there is
the commutator~$a\in P^\prime \setminus Z(G)$ of order~2.
Let $P^\prime \le Z(G)$ and $a\in X_2(S)\setminus Z(G)$, $a^2\in Z(G)$.
Since $G/Z(G)$ is a nonabelian simple group, there is~\cite[3.8.2]{Gor}
an element $g\in G$ such that $\langle aZ(G),a^gZ(G)\rangle$ is not a 2-group.
But a subgroup generated by two involutions is a dihedral group.
So there is an element $tZ(G)\in \langle aZ(G),a^gZ(G)\rangle$ of
odd  prime order~$r$ such that $\langle aZ(G),tZ(G)\rangle =D/Z(G)$
is a dihedral group of order $2r$.
Let $R$ be a Sylow $r$-subgroup of $D$. Since
$RZ(G)=R\times Z(G)_{r^\prime}$ is normal in~$D$,
it follows that $R$ is normal in $D$. In view of Lemma~4,
$a\in C_G(R)$ and $D/Z(G)$ is abelian. This is a contradiction,
and $G$ is soluble.

By induction, all proper subgroups of $G$ have the nilpotent derived subgroups,
in particular, $G^\prime$ is not nilpotent, but $G^{\prime \prime}$ is nilpotent.
Suppose that $S$ is a Schmidt subgroup in $G^\prime$,
$S=P\leftthreetimes Q$, where $P$ is a normal in~$S$ Sylow $p$-subgroup,
$Q$ is a nonnormal in $S$ cyclic Sylow $q$-subgroup, ${S^\prime =P}$.
It is clear that $P\le G^{\prime \prime}$ and $P$ is subnormal in~$G$.
This implies that $P\le O_p(G)$. Here $O_p(G)$ is the largest normal $p$-subgroup
of~$G$. Let $Q_1$ be a Sylow $q$-subgroup of~$G$ containing $Q$.
In view of the Focal Subgroup Theorem \cite[7.3.4]{Gor},
$Q_1\cap G^\prime$ is generated by elements
$x^{-1}y$, where $x,y\in Q_1$ and $x$ conjugates to $y$ in~$G$.
From $(1)$ it follows that
$Q_1\cap G^\prime$ is generated by elements of~$Q_1\cap X_q(G)$.
If $x\in Q_1\cap X_q(G)$, then $x\in C_G(O_p(G))\le C_G(P)$ by Lemma~3, and
$Q\le Q_1\cap G^\prime \le C_G(P)$. This is a contradiction
since $S$ is nilpotent. $\hfill  \square$

\medskip

{\bf Corollary~2.1.} {\sl The derived subgroup of a finite group
is nilpotent if and only if all primary commutators of coprime orders
are permutable.} $\hfill  \square$

\medskip

The following nilpotency criterion for a finite group was obtained
by Baumslag and Wiegold~\cite{BW}.

\medskip

{\bf Theorem 4.} {\sl
Let $G$ be a finite group in which $|ab|=|a||b|$ for elements
$a$, $b$ of coprime orders. Then $G$ is nilpotent.}

\medskip

Here we can also limit to only primary elements.
An element $a$ of a group $G$ is called a $p$\nobreakdash-\hspace{0pt}element
if $|a|=p^n$ for  a non-negative integer $n$.
If the order of $a$ is not divided by~$p$,
then $a$ is a $p^\prime$\nobreakdash-\hspace{0pt}element.
A group $G$ is said to be $p$\nobreakdash-\hspace{0pt}decomposable
if $G$ has a normal Sylow $p$-subgroup and $G$ is
$p$\nobreakdash-\hspace{0pt}nilpotent.

\medskip

{\bf Theorem 5.} {\sl
Let $G$ be a finite group. Fix a prime $p\in \pi(G)$.
A group~$G$ is $p$-decomposable if and only if $|ab|\ge |a||b|$
for a $p$-element $a$ and a primary $p^\prime$-element $b$.}

\medskip

{\sc Proof.}
If $a$ is a  $p$-element and $b$ is a primary $p^\prime$-element of
a finite $p$\nobreakdash-\hspace{0pt}decomposable group, then $|ab|=|a||b|$.

Conversely, assume that the result is not true and
let $G$ be a counterexample of least order.
It is clear that the hypothesis holds for every subgroup of $G$.
Hence all proper subgroups of $G$ are $p$\nobreakdash-\hspace{0pt}decomposable.
Suppose that $G$ is not $p$\nobreakdash-\hspace{0pt}nilpotent.
In view of Ito Theorem~\cite[IV.5.4]{Hup},
$G=P\leftthreetimes Q$ is a Schmidt group,
where $P$ is a normal in~$G$ Sylow $p$-subgroup,
$Q=\langle y\rangle$ is a Sylow subgroup of order~$q^t$.
If $x$ is an element of $P$, then $x^{-1}yx=ay^s$ for some $a\in P$
and a non-negative integer~$s$. Since
$$
a=x^{-1}yx(y^s)^{-1}=x^{-1}yxy^{-1}y(y^s)^{-1}=[x,y^{-1}]y(y^s)^{-1},
\ \  [x,y^{-1}]\in G^\prime \le P, \eqno (2)
$$
this implies that $y(y^s)^{-1}=1$, $y=y^s$ and $x^{-1}yx=ay$.
By the hypothesis,
$$
|y|=|x^{-1}yx|=|ay|\ge |a||y|, \ \ a=1, \ \ x^{-1}yx=y. \eqno (3)
$$
Consequently, $Q$ is normal in $G$, a contradiction.
Thus, $G=H\leftthreetimes P$ is $p$-nilpotent.

If $P$ is not cyclic, then in~$P$ there are two different maximal
subgroups~$P_1$ and~$P_2$. Subgroups $HP_1$ and $HP_2$ are $p$-decomposable,
it follows that $G$ is $p$-decomposable. Let $P=\langle y\rangle$ be cyclic.
If $x$ is an element of~$H$, then $x^{-1}yx=ay^s$ for some $a\in P$
and a non-negative integer~~$s$. Repeating $(2)$ and $(3)$,
we have that $G$ is $p$-decomposable.
$\hfill \square$

\medskip

\textbf{Corollary 5.1.} {\sl
A finite group is nilpotent if and only if ${|ab|\ge |a||b|}$
for any two primary elements $a$ and $b$ of coprime orders.}~$\hfill \square$

\bigskip

Department of mathematics, 

Gomel F. Scorina State University, 

Gomel, BELARUS

\medskip

{\sl E-mail}: Victor.Monakhov@gmail.com

\end{document}